\documentclass[11pt,centertags,reqno,twoside]{amsart}
\usepackage{amssymb}
\usepackage{mathptmx}
\usepackage[mathcal]{euscript}

\usepackage[english]{babel}

\DeclareSymbolFont{SY}{U}{psy}{m}{n}
\DeclareMathSymbol{\emptyset}{\mathord}{SY}{'306}

\newcommand{\bI}{\textbf{1}}

\newcommand{\dist}{\operatorname{dist}}

\allowdisplaybreaks

\newtheorem{theorem}{Theorem}
\newtheorem{corollary}[theorem]{Corollary}

\theoremstyle{definition}

\theoremstyle{remark}
{\it}{\rm}

\begin{document}
%%%%%%%%%%%%%%%%%%%%%%%%%%%%%%%%%%%%%%%%%%%%%%%%%%%%%%%%%%%%%%%%%

\title[Unconditional bases of subspaces]
{Unconditional bases of subspaces related to non-self-adjoint
perturbations of self-adjoint operators}
\thanks{$^*$This work was
supported by the Deutsche For\-sch\-ungs\-gemeinschaft (DFG) and the
Russian Foundation for Basic Research, grant No. 15-51-12389.}

\author[A. K. Motovilov and A. A. Shkalikov]
{Alexander K. Motovilov and Andrei A. Shkalikov}

\address{Alexander K. Motovilov\newline\hspace*{9mm} Joint Institute for Nuclear Research
and Dubna State University, Dubna, Russia}
\email{motovilv@theor.jinr.ru}

\address{Andrei A. Shkalikov\newline\hspace*{9mm} Faculty of Mathematics and Mechanics,
 Lomonosov Moscow  State University, Moscow, Russia}
\email{shkalikov@mi.ras.ru}

%\date{January 24, 2017}

\keywords{Riesz basis, unconditional basis of subspaces,
non-self-adjoint perturbations}.

%%%%%%%%%%%%%%%%%%%%%%%%%%%%%%%%%%%%%%%%%%%%%%%%%%%%%%%%
\begin{abstract}
Assume that $T$ is a self-adjoint operator on a Hilbert space
$\mathcal{H}$ and that the spectrum of $T$ is confined in the union
$\bigcup_{j\in J}\Delta_j$,\, $J\subseteq\mathbb{Z}$, of segments
$\Delta_j=[\alpha_j, \beta_j]\subset\mathbb{R}$ such that
$\alpha_{j+1}>\beta_j$ and
$$
\inf_{j} \left(\alpha_{j+1}-\beta_j\right) = d > 0.
$$
If $B$ is a bounded (in general non-self-adjoint) perturbation of
$T$ with $\|B\|=:b<d/2$ then the spectrum of the perturbed operator
$A=T+B$ lies in the union $\bigcup_{j\in J} U_{b}(\Delta_j)$ of the
mutually disjoint closed $b$-neighborhoods $U_{b}(\Delta_j)$ of the
segments $\Delta_j$ in $\mathbb{C}$. Let $Q_j$ be the Riesz
projection onto the invariant subspace of $A$ corresponding to the
part of the spectrum of $A$ lying in $U_{b}\left(\Delta_j\right)$,
$j\in J$. Our main result is as follows: {\sl  The subspaces
$\mathcal{L}_j=Q_j(\mathcal H)$, $j\in J$, form an unconditional basis
in the whole space $\mathcal H$.}
\end{abstract}
%%%%%%%%%%%%%%%%%%%%%%%%%%%%%%%%%%%%%%%%%%%%%%%%%%%%%

\maketitle

%%%%%%%%%%%%%%%%%%%%%%%%%%%%%%%%%%%%%%%%%%%%%%%%%%%%%
\section{Introduction and main result}
\label{SIntro}
%%%%%%%%%%%%%%%%%%%%%%%%%%%%%%%%%%%%%%%%%%%%%%%%%%%%%%%

We begin with recalling some definitions (see \cite[Ch.6.5]{GK}). A
sequence of nonzero subspaces $\left\{ \mathcal{L}_k \right\}$ of a
Hilbert spaces $\mathcal{H}$ is said to be {\it a basis} if any
element $x\in \mathcal{H}$ is uniquely represented by the series
\begin{equation} \label{xksum}
 x=\sum_k x_k,\quad\text{where \,\,} x_k\in\mathcal{L}_k,
\end{equation}
that converges in the norm of $\mathcal{H}$. A basis of subspaces
$\left\{ \mathcal{L}_k \right\}$ is said to be
\textit{unconditional} if the series \eqref{xksum} converges to $x$
after any rearrangement of its elements. An unconditional basis of
subspaces is also called a Riesz basis of subspaces. In the case
where all the subspaces $\mathcal{L}_j$ are one-dimensional
(finite-dimensional) we can choose elements $y_j\in \mathcal
L_j=Q_j(\mathcal H)$ (a basis $\{y_{k_j}\}$ in $\mathcal L_j$). Then
the sequence of the subspaces $\left\{ \mathcal{L}_k \right\}$ is a
basis if and only if the corresponding system is a  basis (a  basis
with parentheses) in $\mathcal H$.

When one deals with bases of subspaces, it is convenient to work in
terms of projections. We will use definitions and results presented
in \cite[\S 6]{Shka1}. Let $J$ be a finite or infinite ordered set
of indices, $J\subset \mathbb{Z}$. A system of projections
$\{Q_j\}_{j\in J}$ is said to be \textit{complete} if the equalities
$$
\left(Q_j x,y\right)=0, \quad \text{for any} \quad x\in
\mathcal{H}\quad \text{and any}\quad  j\in J,
$$
imply $y=0$. It is easily seen that the system $\{Q_j\}_{j\in J}$ is
complete if and only if any element $x\in\mathcal H$ can be
approximated with an arbitrary accuracy by linear combinations of
elements $x_k \in\mathcal{L}_k =Q_k(\mathcal H)$, $k\in J$.

A system of projections $\{Q_j\}_{j\in J}$ is called
\textit{minimal} if
\begin{equation*}%\label{Qmin}
Q_jQ_k=\delta_{kj}Q_j\quad\text{for any \,\,} j,k \in J.
\end{equation*}
It follows directly from the definitions that if $\{Q_j\}_{j\in J}$
is a minimal system of projections, then the sequence of the
subspaces $\mathcal{L}_j=Q_j(\mathcal{H})$ forms a basis (an
unconditional basis) if and only if the series $\sum_{j\in J}Q_j$
converges (converges after any rearrangement of the indices) in the
strong operator topology to the identity operator.

Further we will make use of the following results (for the
corresponding proofs see \cite[Ch. 6]{GK} and \cite[\S 6]{Shka1}).
\smallskip

\noindent\textbf{Theorem A.} {\sl  Let $\{Q_j\}_{j\in J}$ be a
system of projections in a Hilbert space $\mathcal{H}$. The
following statements are equivalent:
\begin{enumerate}

\item A sequence of subspaces $\mathcal{L}_j=Q_j(\mathcal{H})$, $j \in J$,
is an unconditional basis of the Hilbert space $\mathcal{H}$.

\item There exists an equivalent inner product
in $\mathcal{H}$ such that a sequence of subspaces
$\mathcal{L}_j=Q_j(\mathcal{H})$, $j \in J$, is complete and
mutually orthogonal ($\mathcal{L}_k \perp \mathcal{L}_j$ for $k\ne
j$).

\item There exists a bounded and boundedly invertible operator $K$ in $\mathcal{H}$,
and a complete and minimal system of orthogonal projections $\{P_j\}
$ such that
$$Q_j=K^{-1}P_jK  \qquad j\in J.$$

\item The system of projections $\{Q_j\}_{j\in J}$ is complete, minimal, and the series
$\sum_{j\in J}Q_j$ converges unconditionally.

\item The system of projections $\{Q_j\}_{j\in J}$ is complete, minimal, and
\begin{equation}\label{0}
\sum_{j\in J}\left|\left(Q_jx,x\right)\right|<\infty \qquad
\text{for any \,\,} x \in \mathcal{H}.
\end{equation}

\end{enumerate}
}

Now we are ready to formulate the main result of the paper. In what
follows it is assumed that either the set of indices coincides with
$\mathbb N$ or with $\mathbb Z$.
\medskip

\begin{theorem}
\label{Thm1} {\sl Let $T$ be a self-adjoint operator on a Hilbert
space $\mathcal{H}$. Assume that the spectrum of $T$ is confined in
the union $\Delta:=\bigcup_{j\in\ J}\Delta_j$ of the segments
$\Delta_j=[\alpha_j, \beta_j]\subset\mathbb{R}$ such that
$\alpha_{j+1}>\beta_j$ for all $j\in J$. Assume in addition that
\begin{equation}\label{00}
\inf_{j\in J} \left(\alpha_{j+1}-\beta_j\right) = d > 0.
\end{equation}
Let B be a bounded (generally non-self-adjoint) operator on
$\mathcal{H}$ with $\|B\|=:b<d/2$. Then the spectrum the operator
$A=T+B$ lies in the union $\bigcup_{j\in J} U_{b}(\Delta_j)$ of the
mutually disjoint closed $b$-neighborhoods $U_{b}(\Delta_j)$ of the
segments $\Delta_j$ in $\mathbb{C}$. If $Q_j$, $j\in J$, are the
Riesz projections onto the invariant subspaces $\mathcal{L}_j$ of
$A$ corresponding to the isolated components of its spectrum lying
in $U_{b} \left(\Delta_j\right)$, then the invariant subspaces
$\mathcal{L}_j$, $\mathcal{L}_j=Q_j(\mathcal{H})$, $j\in J$, form an
unconditional basis  in $\mathcal{H}$.}
\end{theorem}
\medskip

There are many papers devoted to the Riesz basis property of the
root vectors of non-self-adjoint operators, which are perturbations
of self-adjoint ones.  The corresponding results and references can
be found in the book of Markus \cite{Ma1} and in the paper of
Shkalikov \cite{Shka1}. First results related to Theorem \ref{Thm1} were
obtained by  Markus \cite{Ma2} and Kato (see \cite[Ch.5, Theorem
4.15a]{Ka}): {\sl  Let $T$ be a self-adjoint operator with discrete
spectrum on a Hilbert space $\mathcal H$  such that its eigenvalues
$\{\lambda_j\}$  are simple  and subject the condition
$\lambda_{j+1}- \lambda_j \to \infty$,  as $j\to \infty$. Then for
any bounded (generally non-self-adjoint) operator $B$ the operator
$T+B$ has discrete spectrum and its root vectors form a Riesz basis
in the space $\mathcal H$. }

A more general result follows from the Markus-Matsaev theorem
\cite[Ch.1, Theorem 6.12]{Ma1}: {\sl Let $T$ be a self-adjoint
operator in a Hilbert space $\mathcal H$, having a finite order
(i.e. its eigenvalues, counting multiplicities, are subject to the
condition $|\lambda_j| \ge C j^p$ with some constants $C, p>0$), and
there exist gaps in the spectrum of $T$ with the lengths $\ge d$.
If $\|B\| \le d/2$ then the the root subspaces of the perturbed
operator $T+B$ form a Riesz basis with parentheses $\mathcal H$ (or
a Riesz basis consisting of finite-dimensional root subspaces).}

The condition for $T$ to be of finite order was dropped in
\cite{Shka1}. However, the compactness of the resolvent
$(T-\lambda)^{-1}$ was essentially used in the proof. To the best of
our knowledge, Theorem \ref{Thm1} is, apparently, the first result in this
topic which deals with an unperturbed operator $T$ possibly having a
non-discrete spectrum.

We have an additional motivation to prove Theorem \ref{Thm1}. We
expect that this result might help to resolve some open problems
concerning bounded perturbations of self-adjoint operators in spaces
with indefinite metric (see \cite{AMSh,AMT}).

\section{Proof of Theorem \ref{Thm1}}
We divide the proof into two parts that are called below Step 1 and
Step 2, respectively. At Step 1 we prove that the sequence of
subspaces $\{\mathcal{L}_j\}_{j\in J}$ forms a basis of
$\mathcal{H}$. At Step 2 we prove that it is, actually, an
unconditional basis.
\smallskip

\noindent
%%%%{\bf Proof.}
\textit{Step 1.} For the sake of definiteness, we assume that the index
set $J$ coincides with the set of all entire numbers
$\mathbb{Z}=\{\dots, -2, -1,0, 1, 2, \dots\}$. This corresponds to
the most general case when infinitely many segments $\Delta_j$ lie
both on $\mathbb{R}^-$ and $\mathbb{R}^+$.

Under the hypothesis $\|B\|=b$ with $0\leq b<d/2$ that we assume,
the closed neighborhoods $U_{b}(\Delta_j)$ of the segments
$\Delta_j$, $j\in J$, are disjoint and, surely, $U_{b}(\Delta)
=\bigcup_{j\in J} U_{b}(\Delta_j)$. The first assertion of the
theorem on the inclusion of the spectrum $\sigma(A)$ of $A$ in the
union $\bigcup_{j\in J} U_{b}(\Delta_j)$ is a consequence of the
well-known estimate
\begin{equation}
\label{eq0} \left\| \left( T-\lambda\right) ^{-1} \right\|  \leq
\frac{1}{\dist (\lambda, \sigma(T))}
 \leq \frac{1}{\dist (\lambda,\Delta)}, \qquad \Delta: = \bigcup_{j\in J} \Delta_j,
\end{equation}
where $\sigma(T)$, $\sigma(T)\subset\Delta$, is the spectrum of the
self-adjoint operator $T$. For any $\lambda$ lying outside
$U_{b}\left(\Delta\right)$ we  have
\begin{equation}
\label{eq1} \delta:=\dist (\lambda,\Delta)>b.
\end{equation}
Then combining \eqref{eq0} and \eqref{eq1} with the bound
\begin{equation}\label{S}
\bigl\|\bigl((\bI+B(T-\lambda)^{-1}\bigr)^{-1}\bigr\|\leq
\frac{1}{1-b/\delta}
\end{equation}
for the resolvent $(A-\lambda)^{-1}=(T+B-\lambda)^{-1}$ one finds
\begin{align}
\nonumber \|(A-\lambda)^{-1}\|=&\left\|\left(T-\lambda
\right)^{-1} \left(\bI+B(T-\lambda)^{-1}\right)^{-1} \right\|\\
\label{AdelB} &\leq \frac{1}{\delta-b} < \infty,
\end{align}
where the quantity $\delta=\delta(\lambda)$ is defined in
\eqref{eq1}. Hence, any $\lambda\in\mathbb{C}\setminus
U_{b}(\Delta)$ belongs to the resolvent set of the perturbed
operator $A=T+B$ and then $\sigma(A)\subset
U_{b}(\Delta)=\bigcup_{j\in J} U_{b}(\Delta_j)$. Since the
neighborhoods $U_{b}(\Delta_j)$ for different $j\in J$ are disjoint,
namely,
$$
\dist(U_{b}(\Delta_j),U_{b}(\Delta_k))\geq d-2 b,\quad j\neq k,
$$
the spectral sets of $A$ confined in $U_{b}(\Delta_j)$ are also
disjoint. The Riesz projections $Q_j$ for these sets are well
defined. In particular, given an arbitrary $b'\in (b, d/2)$, one
may write $Q_j$ in the form
\begin{equation}\label{1}
Q_j=-\frac{1}{2\pi i} \int_{\Gamma_j}
\left(A-\lambda\right)^{-1}\,d\lambda, \qquad \Gamma_j=\partial
U_{b'} \left(\Delta_j\right),
\end{equation}
where the integration along contour $\Gamma_j$ is performed in the
anti-clockwise direction.

Below we will also use the representation
\begin{equation}\label{2}
\left(A-\lambda \right)^{-1}=\left(T-\lambda
\right)^{-1}-G(\lambda), \quad \lambda\notin
U_{b}\left(\Delta\right),
\end{equation}
where
\begin{equation}\label{3}
G(\lambda)= (A-\lambda)^{-1}B\left(T-\lambda \right)^{-1}=
\left(T-\lambda \right)^{-1} M(\lambda) B\left(T-\lambda
\right)^{-1}
\end{equation}
and
$$
M(\lambda) =  \left(\bI+B(T-\lambda)^{-1}\right)^{-1}.
$$
Now denote by $R_n$ the rectangle in $\mathbb{C}$ whose vertical
sides pass through the points
\begin{equation}\label{c_n}
c_{-n}=(\beta_{-n-1}+\alpha_{-n})/2 \ \, \text{and} \
c_{n}=(\beta_{n}+\alpha_{n+1})/2
\end{equation}
 while the horizontal sides
coincide with the segments $[c_{-n}\pm i\gamma_n, c_{n}\pm
i\gamma_n]$ where
\begin{equation}
\label{gamn} \gamma_n=\max\{|c_{-n}|,|c_n|\}.
\end{equation}
Clearly, the sides of $R_n$ do not intersect the set
$U_{b}\left(\Delta\right)$. By virtue of \eqref{2} one obtains
\begin{equation}\label{Q}
\sum_{j=-n}^n Q_jx=-\frac{1}{2\pi i}\int_{\partial
R_n}\left(A-\lambda \right)^{-1}x\,\, d\lambda=\sum_{-n}^nP_jx+I_n
x,
\end{equation}
where $I_n$ are the respective contour integrals of the
operator-valued function $G(\lambda)$ along $\partial R_n$,
\begin{equation}
\label{In} I_n:=\frac{1}{2\pi i}\int_{\partial
R_n}G(\lambda)\,d\lambda,
\end{equation}
and $P_j$ are the spectral projections onto the spectral subspaces
of the self-adjoint operator $T$ associated with the parts of its
spectrum inside the corresponding segments $\Delta_j$.

Given an arbitrary $x\in\mathcal{H}$ we have
\begin{equation}
\label{sump} \sum\limits_{j=-n}^n P_jx\to x\quad\text{as \,\,} n\to
\infty.
\end{equation}
Thus, in order to prove that $\sum\limits_{j=-n}^n Q_jx\to x$ one
only needs to show that the sequence of $I_n x$ in \eqref{Q}
converges to zero as $n\to\infty$.

First, let us show that the operators $I_n$ are uniformly bounded.
We have $I_n = I_n^1+I_n^2$, where $I_n^1$ and $I_n^2$  are the
integrals along the horizontal and vertical sides of the rectangles
$R_n$, respectively. For $\lambda$ varying on the horizontal  sides
of the rectangle $\partial R_n$ we have the estimate
\begin{equation*}
 \|(T-\lambda)^{-1}\|\leq\frac{1}{\gamma_n}, \quad
\lambda=\xi\pm i\gamma_n,\quad \xi\in[c_{-n},c_n],
\end{equation*}
where $\gamma_n$ is defined by \eqref{gamn}. Hence, by virtue of
\eqref{AdelB} and  \eqref{3}, it follows
\begin{equation}\label{Bhori}
 \|G(\lambda)\|\leq\frac{b}{(\gamma_n -b)\gamma_n}, \quad
\lambda=\xi\pm i\gamma_n,\quad \xi\in[c_{-n},c_n].
\end{equation}
Taking into account that the lengths of the horizontal sides do not
exceed $2\gamma_n$ and  $\gamma_n \to \infty$ as $n\to \infty$,  we
get $\|I^1_n\| \to 0$ as $n\to\infty$.

Let us estimate the norms of the operators $I^2_n$. For the sake of
definiteness consider the right vertical side of the rectangle $R_n$
and divide it in three parts
$$
\omega_n \cup \omega^+_n \cup \omega^-_n, \ \ \text{where} \
\omega_n = (c_n-id, c_n+id),\ \ \omega^\pm = [c_n\pm id,
c_n\pm i\gamma_n].
$$
The lengths of the intervals $\omega_n$ equal $2d$ and due to
\eqref{AdelB} and \eqref{3}
\begin{equation}\label{Om1}
\left\|\int_{\omega_n} G(\lambda)\, d\lambda\right\| \leq
\frac{2bd}{(d/2-b)(d/2)}.
\end{equation}
For $\lambda\in \omega^\pm_n$ we have
\begin{equation}
\label{dbIm}
\delta(\lambda)-b \geq |\mathop{\rm Im}\lambda|-b.
\end{equation}
Thus, again applying \eqref{AdelB} and \eqref{3}, one obtains
\begin{equation}\label{Om2}
\left\|\int_{\omega^\pm_n} G(\lambda)\, d\lambda\right\| \leq
 \int_{d}^{\gamma_n}\, \frac{b}{\tau(\tau-b)}\, d\tau \leq
 \int_{d}^{\gamma_n}\, \frac{b}{(\tau-b)^2}\, d\tau\leq \frac b{d-b}< 1.
\end{equation}
Therefore,
\begin{equation}\label{CIn}
\|I_n\| \leq C,
\end{equation}
where the constant $C$ depends only on $d$ and $b$.

Let us show that $\|I_nx\| \to 0$ as $n\to \infty$ for any fixed
$x\in \mathcal H$. To this end, choose  some $\varepsilon>0$ and,
first, find $N\in\mathbb{N}$ such that
\begin{equation}
\label{xNeps} \|x-x_N\|<\frac{\varepsilon}{2C},
\end{equation}
where
$$
x_N=\sum_{j=-N}^{N}  P_j x.
$$
Obviously, from \eqref{CIn} and \eqref{xNeps} it follows
\begin{equation}
\label{IxN} \|I_n(x-x_N)\|<\frac{\varepsilon}{2}\quad\text{for any
\,}n\in\mathbb{N}.
\end{equation}
Let us estimate  $\|I_n x_N\|$ as $n\to \infty$. Denote
$e(t):=\left(E(t)x,x\right)$ where  $E(t)$ is the spectral function
of the self-adjoint operator $T$, and observe that, by the spectral
theorem,
\begin{equation}
\label{Tlb} \left\|\left(T-\lambda \right)^{-1}x_N\right\|^2=
\int\limits_{c_{-N}}^{c_N}\frac{de(t)}{(t-\xi)^2+\tau^2}, \quad
\lambda=\xi+i\tau \,
 \not\in\sigma(T).
\end{equation}
For $N$ fixed the equality \eqref{Tlb} implies
\begin{equation}
\label{Tlam} \left\|\left(T-\lambda
\right)^{-1}x_N\right\|=O\left(|\lambda|^{-1}\right) \ \ \text{ as}\
\, |\lambda|\to \infty.
\end{equation}
Now we can modify estimate \eqref{Om1}  and get from \eqref{AdelB}
and \eqref{Tlam}
$$
\left\|\int_{\omega_n} G(\lambda)x_N\, d\lambda\right\|
=O\left(|c_n|^{-1}\right)\to 0\ \  \text{as} \ n\to \infty.
$$
Analogously, by taking into account \eqref{dbIm} we can modify the
estimate \eqref{Om2} and obtain
\begin{equation*}
\left\|\int_{\omega^\pm_n} G(\lambda)x_N\, d\lambda\right\| \leq
 \int_{d}^{\gamma_n}\, \frac{d\tau}{(\tau-b)\, \sqrt{(c_n-c_N)^2 +\tau^2}}
 = O\left(\frac{\ln c_n}{c_n}\right) = o(1) \ \  \text{as} \  n\to \infty.
\end{equation*}
 The last two estimates together give
\begin{equation*}
I_n x_N\to 0\quad \text{as \,} n\to\infty.
\end{equation*}
Hence, there is $N_1\in\mathbb{N}$ $(N_1\geq N)$ such that $\|I_n
x_N\|<\varepsilon/2$ whenever $n>N_1$. Taken together with
\eqref{IxN} this yields $\|I_n x\|<\varepsilon$ for  $n>N_1$.
Therefore, we have proven that the sequence of $I_n$ strongly
converges to zero as $n\to\infty$ and then from \eqref{Q} and
\eqref{sump} it follows that for any $x\in\mathcal{H}$ the two-sided
series $\sum_{j=-n}^n Q_j x$ converges to $x$. Thus, the system of
subspaces $\mathcal{L}_j=Q_j(\mathcal{H})$, $j\in\mathbb{Z}$ is
complete. That these subspaces are linearly independent follows from
the mutual orthogonality of the Riez projections \eqref{1} in the
sense that $Q_j Q_k=\delta_{jk}Q_j$ for any $j,k\in\mathbb{Z}$ (see,
e.g. \cite[Ch. I \S 1.3]{GK}). Thus, the system
$\{\mathcal{L}_j\}_{j=-\infty}^\infty$ represents a basis of
subspaces in $\mathcal{H}$.

\medskip

\noindent\textit{Step 2.} By Theorem A, in order to prove that the
above basis of subspaces $\{\mathcal{L}_j\}_{j=-\infty}^\infty$ is
unconditional,  it suffices to show that the series
$\sum_{j=-\infty}^\infty\left|\left(Q_jx,x\right)\right|$ converges
for any $x \in \mathcal{H}$. First, let us transform the contours
$\Gamma_j$ in the integrals \eqref{1} into the contours
$\widetilde{\Gamma}_j$, which surround the rectangles with the
vertical sides $\lambda=c_{j-1}+i\tau$ and $\lambda =c_j +i\tau$ where
$\tau$ is varying in $[-d,d]$, and the horizontal sides $\xi\pm id$,
$\xi\in[c_{j-1}, c_j]$. As previously, the numbers $c_{j}$ are
defined by \eqref{c_n} and coincide with the centers
of the gaps between the neighboring intervals
$\Delta_j$ and $\Delta_{j+1}$.

Taking  into account the representation \eqref{3} and  estimate
\eqref{S}, we immediately conclude that
$$
2\pi \sum_{j=-\infty}^{\infty}|(Q_jx,x)|\leq
\sum_{j=-\infty}^{\infty}\left| \int_{\widetilde{\Gamma}_j}\left(
\left(T-\lambda \right)^{-1}x,x\right)\,
d\lambda\right|+\sum_{j=-\infty}^{\infty}\left|
\int_{\widetilde{\Gamma}_j}\left(G(\lambda)x,x\right)\,
d\lambda\right|.
$$
The first series converges since it simply coincides with the series
$\sum_{j\in\mathbb{Z}}(P_jx,x)=\|x\|^2$; we recall that $P_j$ are
the spectral projections of the self-adjoint operator $T$ associated
with the segments $\Delta_j$. By virtue of \eqref{S} and \eqref{3},
the second series can be estimated as follows:
\begin{align}
\nonumber
\sum_{j=-\infty}^\infty\left|\int_{\widetilde{\Gamma}_j}\left(G(\lambda)x,x\right)\,
d\lambda\right|&\leq
C_1\sum_{j=-\infty}^\infty\int_{\widetilde{\Gamma}_j}\left\|(T-\lambda)^{-1}x\right\|^2
\left|d\lambda\right|\\
\label{6} &\leq
C_1\left(\int_{\Gamma_+\cup\Gamma_-}\left\|\left(T-\lambda
\right)^{-1}x\right\|^2 \left|
d\lambda\right|+2\sum_{j\in\mathbb{Z}}\int_{\omega_j}\left\|
\left(T-\lambda \right)^{-1}x\right\|^2
\left|d\lambda\right|\right),
\end{align}
where $\Gamma_{\pm}$ are the lines $\lambda=\xi\pm id$,\,
$\xi\in\mathbb{R}$, $\omega_j$ are the vertical segments
$\lambda=c_j+i\tau$, $-d\leq \tau\leq d$, and $C_1= const $.
% (cf. \eqref{Sbound}).
Convergence of the integrals and the series in \eqref{6} can be
proven by using the spectral theorem. As before, denote
$e(t):=\left(E(t)x,x\right)$ where $E(t)$ stands for the spectral
function of $T$. Then
\begin{align}
\nonumber \int_{\Gamma_\pm}\left\|\left(T-\lambda
\right)^{-1}x\right\|^2 \left|d\lambda\right|&= \nonumber
\int_{\Gamma_\pm}|d\lambda|\int_{\mathbb{R}} \frac{de(t)}{|t-\lambda|^2}\\
\nonumber
&= \int_{\mathbb R} de(t)\int_{\mathbb{R}} \frac{d\xi}{|t-\xi|^2+d^2}\\
\nonumber
&= \frac{\pi}{d} \int_{\mathbb{R}} de(t)\\
\label{7} &=\frac{\pi}{d}\|x\|^2.
\end{align}
Further, notice that for $c_j$ given by \eqref{c_n}, the lower bound
\eqref{00} yields
$$
|t-c_j|\geq \ \,  \left\{\begin{matrix} d/2+d(j-k) \quad \text{for
all}\quad  t\in \Delta_k,\ \,  k\leq
j,\\
d/2+d(k-j-1) \quad \text{for all}\quad  t\in \Delta_k,\ \,  k\geq
j+1.
\end{matrix}\right.$$
Hence, for $\lambda = c_j+i\tau$, $\tau\in[-d,d]$, we have
\begin{align}
\nonumber
\sum_{j\in\mathbb{Z}}\int_{\gamma_\pm}\left\|\left(T-\lambda\right)^{-1}x\right\|^2
\left|d\lambda\right|&= \sum_{j\in\mathbb{Z}}
\int\limits_{-d}^d d\tau \int_{\mathbb{R}} \frac{de(t)}{|t-c_j-i\tau|^2}\\
\nonumber &=\sum_{j\in\mathbb{Z}} \int\limits_{-d}^d d\tau
\sum_{k\in\mathbb{Z}} \int_{\Delta_k}
\frac{de(t)}{|t-c_j-i\tau|^2} \\
\nonumber &\quad\leq  2d\sum_{j\in\mathbb{Z}}\sum_{k\in\mathbb{Z}}
\int_{\Delta_k}  \frac{de(t)}{|t-c_j|^2}\\
\nonumber &\quad\leq\sum_{k\in\mathbb{Z}} \|P_k
x\|^2\left(2\left(\frac{2}{d}\right)^2+\sum_{j\in\mathbb{Z},\,\,
j\neq k}
\frac{1}{d^2|k-j|^2}\right)\\
\label{Tlx} &\quad\leq \frac{2C_2}{d^2}\sum_{k\in\mathbb{Z}} \|P_k
x\|^{2}=\frac{2C_2}{d^2}\|x\|^{2},
\end{align}
where
$C_2=4+\sum\limits_{j=1}^\infty\frac{1}{j^2}=4+\frac{\pi^2}{6}$. The
estimate \eqref{6} together with \eqref{7} and \eqref{Tlx} entails
the bound \eqref{0}. This completes the proof of Theorem \ref{Thm1}.
\medskip

The following statement is a simple corollary of Theorem \ref{Thm1} (combined
with Theorem A).

\begin{corollary}
{\sl Assume the hypothesis of Theorem \ref{Thm1}. Then there exists an inner
product $\langle\ , \ \rangle$ in $\mathcal{H}$ with the following
properties.

\begin{enumerate}

\item[(1)] The product $\langle\ , \ \rangle$ is norm-equivalent to the original
inner product $(\ , \ )$ in $\mathcal{H}$.

\item[(2)] The subspaces
$\mathcal{L}_j=Q_j\mathcal{H}$ are mutually orthogonal with respect
to the inner product $\langle\ , \ \rangle$ and, with respect to
$\langle\ , \ \rangle$, the Hilbert space $\mathcal{H}$ admits the
orthogonal decomposition
\begin{equation}
\label{Hdec} \mathcal{H}=\bigoplus\limits_{j\in J}\mathcal{L}_j.
\end{equation}

\item[(3)] The subspaces $\mathcal{L}_j$, $j\in J$, are reducing for
the perturbed operator $A=T+B$ and, with respect to the
decomposition \eqref{Hdec}, this operator admits a block diagonal
matrix representation
$$
A=\mathop{\rm diag}(\ldots,A_{-2}A_{-1},A_0,A_1,A_2,\ldots),
$$
where $A_j=A\big|_{\mathcal{L}_j}$, $j\in J$, denotes the part of
$A$ in the reducing subspace $\mathcal L_j$.
\end{enumerate}
}
\end{corollary}

\end{document}